\documentclass[12pt]{amsart}
\usepackage{amsthm}
\usepackage{amsmath}
\usepackage{amsfonts}
\usepackage{amssymb}
\usepackage{mathrsfs}
\usepackage{mathtools}
\usepackage{color}

\theoremstyle{plain}
\newtheorem{lem}{Lemma}[section]
\newtheorem{prop}[lem]{Proposition}
\newtheorem{thm}[lem]{Theorem}
\newtheorem{cor}[lem]{Corollary}
\theoremstyle{definition}
\newtheorem{defn}[lem]{Definition}

\newtheorem{rem}[lem]{Remark}
\newcommand{\nequation}{\setcounter{equation}{0}}
\renewcommand{\theequation}{\mbox{\arabic{section}.\arabic{equation}}}
\newcommand{\R}{\mathbb R}

\newcommand{\Z}{\mathbb Z}
\newcommand{\N}{\mathbb N}
\newcommand{\Diff}{\text{\rm Diff}}

\newcommand{\id}{\text{\rm id}}

\newcommand{\Ad}{\mbox{\rm Ad}}

\renewcommand{\S}{S^1}

\newcommand{\eps}{\varepsilon}
\newcommand{\ska}[2]{\left\langle #1,#2\right\rangle}

\newcommand{\bea}{\begin{eqnarray}}
\newcommand{\eea}{\end{eqnarray}}
\newcommand{\red}{\textcolor{black}}
\newcommand{\magenta}{\textcolor{black}}
\title[Two-component CH and DP equations]{The geometry of the two-component Camassa--Holm and Degasperis--Procesi equations}
\author{J. Escher, M. Kohlmann, and J. Lenells}
\address{Institute for Applied Mathematics, University of Hannover, D-30167 Hannover, Germany}
\email{escher@ifam.uni-hannover.de}
\email{kohlmann@ifam.uni-hannover.de}
\address{Department of Mathematics, Baylor University, One Bear Place 97328, Waco, USA}
\email{Jonatan\_Lenells@baylor.edu}
\date{\today}
\keywords{Camassa--Holm equation, Degasperis--Procesi equation,
semidirect product, geodesic flow, sectional curvature}
\subjclass[2000]{}

\begin{document}

\begin{abstract}
We use geometric methods to study two natural two-component
generalizations of the periodic Camassa--Holm and\linebreak
Degasperis--Procesi equations. We show that these generalizations
can be regarded as geodesic equations on the semidirect product of
the diffeomorphism group of the circle $\Diff(\S)$ with some space
of sufficiently smooth functions on the circle. Our goals are to
understand the geometric properties of these two-component systems
and to prove local well-posedness in various function spaces.
Furthermore, we perform some explicit curvature calculations for
the two-component Camassa--Holm equation, giving explicit examples
of large subspaces of positive curvature.
\end{abstract}
%
%
%%%%%%%%%%%%%%%%%%%%%%%%%%%%%%%%%%%%%%%%%%%%%%%%%%%%%%%%%%%%%%%%%%%%%%%%%%%%%%%%%%%%%%%%%%%%%%%%%%%%%%%%%%%%%%%%%%%%%%%%%%%
%%%%%%%%%%%%%%%%%%%%%%%%%%%%%%%%%%%%%%%%%%%%%%%%%%%%%%%%%%%%%%%%%%%%%%%%%%%%%%%%%%%%%%%%%%%%%%%%%%%%%%%%%%%%%%%%%%%%%%%%%%%
%%%%%%%%%%%%%%%%%%%%%%%%%%%%%%%%%%%%%%%%%%%%%%%%%%%%%%%%%%%%%%%%%%%%%%%%%%%%%%%%%%%%%%%%%%%%%%%%%%%%%%%%%%%%%%%%%%%%%%%%%%%
%
%
\maketitle

\tableofcontents
\section{Introduction}\label{Introduction}\nequation
In a seminal paper \cite{A66}, Arnold pointed out that the Euler
equations for the motion of a rotating rigid body and the Euler
equations of hydrodynamics can both be viewed geometrically as
geodesic equations on a Lie group endowed with an invariant
metric. More recently, several other equations of physical
interest have been found to arise in a similar way; examples
include the Korteweg--de Vries, Burgers, Camassa--Holm (CH), and other
\red{Euler--Poincar\'e equations}. This geometric viewpoint is not only
aesthetically appealing, but is also useful in the study of
well-posedness and stability issues. It is therefore of interest
to find and study further examples of this type.

The CH equation is a re-expression of the geodesic flow on the
diffeomorphism group of the circle $\Diff(\S)$ equipped with the
$H^1$ right-invariant metric \cite{Mis98, Shk98}. Recently, it has
been demonstrated \cite{EK09} that the Degasperis--Procesi (DP)
equation \cite{DP} also can be recast as a geodesic equation on
$\Diff(S^1)$, although in this case the connection does not derive
from an invariant metric \cite{K09}. Just like CH, the DP equation
is an approximation to the governing equations of motion for the
classical water wave problem in the shallow-water regime cf.
\cite{J02}. Both the CH and DP equations are integrable and admit
peakon solutions \cite{CH, DHH}. The integrability manifests
itself in the existence of a Lax pair and a bi-Hamiltonian
structure for each of the equations.

In this paper, we will develop the geometric picture for the
following two-component generalizations of the CH and DP
equations:
\begin{align}\label{2CH}\tag{2CH}
\left\{\begin{array}{lcl}
  m_t & = & -um_x-2mu_x-\rho\rho_x,\\
    \rho_t & = & -(\rho u)_x, \\
\end{array}\right.
\end{align}
and
\begin{align}\label{2DP}\tag{2DP}
\left\{\begin{array}{cll}
    m_t & = & -3mu_x-m_xu-\rho u_x+2\rho\rho_x, \\
    \rho_t & = & -2\rho u_x -\rho_xu,
  \end{array}\right.
\end{align}
where $u(x,t)$Ê and $\rho(x,t)$ are real-valued functions of $x
\in S^1 \simeq \R/\Z$ and $t \in \R$, and $m = u - u_{xx}$.

The system (\ref{2CH}) was first derived in \cite{OR2006} using
bi-Hamiltonian methods. The system admits a Lax pair formulation
and is integrable. In fact, it is related to the first negative
flow of the AKNS hierarchy via a reciprocal transformation
\cite{CLZ06, F06}. A derivation of (\ref{2CH}) in the context of
shallow water waves appears in \cite{CI2008}. Well-posedness and
blow-up results are obtained in \cite{ELY07, GY2010}.\footnote{In
some of these references the term $-\rho\rho_x$ in (\ref{2CH}) is
chosen to have the opposite sign.}

The system (\ref{2DP}) was first proposed in \cite{Pop06} as a
natural generalization of the DP equation in the context of
supersymmetry. Although the approach of \cite{Pop06} automatically
yields one Hamiltonian structure for (\ref{2DP}), neither a second
Hamiltonian structure nor a Lax pair could be found. The question
of the integrability of (\ref{2DP}) therefore remains open.

We will show that the two-component generalizations \eqref{2CH} and \eqref{2DP} can be regarded
as geodesic equations on the semidirect product $\Diff(\S)
\circledS \mathcal{F}(S^1)$, where $\mathcal{F}(S^1)$ denotes a
space of sufficiently smooth real-valued functions on the circle.
For 2CH, the geodesic equation derives from a natural
right-invariant Riemannian metric, whereas for 2DP the affine
connection is not compatible with any such metric. The geometric
construction will give immediate proofs of local well-posedness
for both systems in $H^{s}(S^1) \times H^{s-1}(S^1)$ or
$C^{n}(S^1) \times C^{n-1}(S^1)$ for sufficiently smooth initial
data. Moreover, we will show that the local well-posedness can be
extended to the Fr\'echet space $C^\infty(S^1) \times
C^\infty(S^1)$. \red{Our main result reads as follows.
\begin{thm}\label{thmmain}
There exist open intervals $J_1$ and $J_2$ centered at $0$ and an open
neighborhood $U$ of $(0,0) \in C^\infty(\S)\times C^\infty(\S)$
such that for each $(u_0,\rho_0)\in U$  there exist a unique
solution
$$(u, \rho) \in C^\infty\bigl(J_1, C^\infty(S^1) \times C^\infty(S^1)\bigr)$$
of (\ref{2CH}) and a unique solution
$$(v, \eta) \in C^\infty\bigl(J_2, C^\infty(S^1) \times C^\infty(S^1)\bigr)$$
of (\ref{2DP}) with $(u(0), \rho(0)) = (v(0),\eta(0)) = (u_0,\rho_0)$. Furthermore,
the solutions depend smoothly on the initial data in the sense
that the local flows
$$\Phi_i\colon J_i \times U \to C^\infty(\S)\times
C^\infty(\S),$$
for $i=1,2$, defined by $\Phi_1(t, u_0, \rho_0) = (u(t; u_0,
\rho_0), \rho(t; u_0, \rho_0))$ and $\Phi_2(t, u_0, \rho_0) = (v(t; u_0,
\rho_0), \eta(t; u_0, \rho_0))$ are smooth maps.
\end{thm}
Although a geometric reformulation of the 2CH as a geodesic flow is presented in \cite{GO06}, see also \cite{G00},
our work contains the following novel aspects: We apply the geometric picture to obtain local well-posedness results (in particular in the smooth category) and we provide a detailed discussion of the sectional curvature associated with the 2CH. A generalization of our approach for 2DP has previously not been presented in the literature.
}

Our paper is organized as follows: In Section \ref{Secsemidirect},
we introduce the relevant function spaces and semidirect products.
In Section~\ref{Sec2CH}, we establish the geometric interpretation
of 2CH as a geodesic equation with respect to a right-invariant
metric and prove local well-posedness in various settings. The 2DP
equation is considered in Section~\ref{Sec2DP}. In
Section~\ref{SecCurvature}, we present some explicit computations
of the sectional curvature for the 2CH equation. In an appendix,
the geometric interpretations of 2CH, CH, and the rotating rigid
body are compared in an attempt to emphasize the unifying features
of the approach.

\section{Function spaces and semidirect products}\label{Secsemidirect} \nequation
We will show that 2CH and 2DP are geodesic equations on the semidirect product
\begin{equation}\label{Gdef}
  G = \Diff(S^1) \circledS \mathcal{F}(S^1),
\end{equation}
where $\Diff(S^1)$ denotes the group of orientation-preserving
diffeomorphisms of the circle $S^1 \simeq \R/\Z$ and
$\mathcal{F}(S^1)$ denotes a space of sufficiently regular
real-valued functions on $S^1$ (the exact regularity assumptions
will be made precise below).

Let $(\varphi, f)$ and $(\psi, g)$ be two elements of $G$. The group product in $G$ is defined by
$$(\varphi, f)(\psi, g)\coloneqq(\varphi\circ\psi, g + f \circ \psi)$$
where $\circ$ denotes composition. The neutral element of $G$ is
$(\id,0)$ and $(\varphi,f)$ has the inverse
$(\varphi^{-1},-f\circ\varphi^{-1})$. Of particular interest to us
will be the right translation operator $R_{(\psi, g)}:G \to G$
defined by
$$R_{(\psi, g)}(\varphi, f) = (\varphi, f)(\psi, g).$$

Several different regularity assumptions can be imposed on the
elements of $G$. The structure of equations (\ref{2CH}) and
(\ref{2DP}) suggests that the function $\rho$ should be allowed to
have one spatial derivative less than $u$. This suggests the
following choice for $G$:
\begin{equation}\label{HsGdef}
  H^sG := H^{s}\Diff(S^1) \circledS H^{s -1}(S^1),
\end{equation}
where $H^{s}\Diff(S^1)$ denotes the space of
orientation-preserving diffeomorphisms of $S^1$ of Sobolev class
$H^s$. We will assume that $s > 5/2$. In this case,
$H^{s}\Diff(S^1)$ is a Hilbert manifold and a topological group
and the composition map
$$(\varphi, f) \mapsto f \circ \varphi: H^{s}\Diff(S^1) \times H^{s-1}(S^1) \to H^{s-1}(S^1)$$
is continuous cf. \cite{EM70}. Thus, $H^sG$ is a topological group
and a smooth manifold modeled on the Hilbert space $H^{s}(S^1)
\times H^{s-1}(S^1)$.

Another natural choice for $G$ is
\begin{equation}\label{CnGdef}
  C^nG := C^{n}\Diff(S^1) \circledS C^{n-1}(S^1),
\end{equation}
where $C^{n}\Diff(S^1)$ denotes the space of
orientation-preserving diffeomorphisms of $S^1$ of class $C^n$. We
will assume that $n \geq 2$. In this case, $C^nG$ is a topological
group and a smooth manifold modeled on the Banach space
$C^{n}(S^1) \times C^{n-1}(S^1)$. Note that $H^sG$ and $C^nG$ are
{\it not} Lie groups, since left multiplication is only continuous
and not smooth.

Finally, we may choose $G$ as
\begin{equation}\label{CinftyGdef}
  C^\infty G := C^\infty\Diff(S^1) \circledS C^\infty(S^1).
\end{equation}
This is a Lie group (the multiplication and inverse maps are
smooth) and a Fr\'echet manifold modeled on $C^\infty(S^1) \times
C^\infty(S^1)$. In contrast to $H^sG$ and $C^nG$, it is {\it not}
a Banach manifold.

The three choices (\ref{HsGdef})-(\ref{CinftyGdef}) for $G$ are
all of interest due to their different advantages. We will first
develop the theory for $H^sG$ and then consider $C^nG$ and
$C^\infty G$.

We refer the reader to \cite{HMR98, HT09} for further information on geodesic flows on semidirect products.

\section{The 2CH equation as a geodesic equation}\label{Sec2CH} \nequation
Let $G$ be the semidirect product defined in (\ref{Gdef}). We will
define a metric $\langle \cdot, \cdot \rangle$ and a compatible
covariant derivative $\nabla$ on $G$ and show that a curve
$(\varphi(t), f(t))$ in $G$ is a geodesic with respect to $\nabla$
if and only if $(u(t), \rho(t)) \in T_{(\id, 0)}G$ defined by
\bea\label{Eulervel} (u,\rho) = TR_{(\varphi, f)^{-1}}(\varphi_t, f_t) = (\varphi_t\circ\varphi^{-1},f_t\circ\varphi^{-1})\eea
satisfies the 2CH equation.

\subsection{The $H^s$-category}
We first consider the $H^s$-setting and let $G$ be the group
$H^sG$, $s > 5/2$, defined in (\ref{HsGdef}). We define a bilinear
operator $\Gamma_{(\text{id},0)}$ on $H^{s}(\S)\times H^{s-1}(\S)$
by
\begin{subequations}\label{Gammadef}
\begin{equation}\label{Christoffel2CH}
\Gamma_{(\text{id},0)}((u,\rho), (v, \tau)) =
\begin{pmatrix}
\Gamma^0_{\text{id}}(u,v)-\frac{1}{2}A^{-1}\partial_x(\rho\tau) \\
  -\frac{1}{2}(u_x\tau + v_x\rho)
\end{pmatrix},
\end{equation}
where $A = 1-\partial_x^2$ and
\begin{equation}
  \Gamma_{\text{id}}^0(u,v) = -A^{-1}\partial_x\left(uv+\frac{1}{2}u_xv_x\right)
\end{equation}
is the Christoffel operator associated with the CH equation (cf. \cite{CK02, L07, Mis98}).
For vector fields $X$ and $Y$ on $H^sG$, we define
\begin{equation}\label{Gammarightinvariance}
\Gamma_{(\varphi,f)}(X,Y)=\Gamma_{(\text{id},0)}(X(\varphi,f)\circ\varphi^{-1},
Y(\varphi,f)\circ\varphi^{-1})\circ\varphi.
\end{equation}
\end{subequations}
Then $\Gamma$ is a right-invariant Christoffel map on $H^sG$\red{, i.e.,
$$TR_{(\psi,g)}[\Gamma_{(\varphi,f)}(X,Y)]=\Gamma_{R_{(\psi,g)}(\varphi,f)}(TR_{(\psi,g)}X(\varphi,f),TR_{(\psi,g)}Y(\varphi,f)),$$
for all $(\varphi,f),(\psi,g)\in H^sG$}. The associated covariant derivative $\nabla$ is defined by
\begin{equation}\label{nabladef}
  (\nabla_XY)(\varphi,f) = DY(\varphi,f)\cdot
X(\varphi,f)-\Gamma_{(\varphi,f)}(Y(\varphi,f),X(\varphi,f)).
\end{equation}
\red{Observe that the Christoffel map $\Gamma$ is the infinite-dimensional analog
of the Christoffel symbols $\Gamma^i_{jk}$ familiar from
finite-dimensional differential geometry (see \cite{Lang2}; our
$\Gamma$ is denoted by $B$ in \cite{Lang2}). Furthermore, it follows immediately from definition (\ref{nabladef}) that
$\nabla$ is a \magenta{torsionless} covariant derivative in the sense that
\begin{enumerate}\renewcommand{\labelenumi}{(\emph{\roman{enumi}})}
\item $\nabla_{fX}Y=f\nabla_XY$,
\item $\nabla_XY-\nabla_YX=[X,Y]$,
\item $\nabla_X(fY)=(Xf)Y+f\nabla_XY$,
\end{enumerate}
for all vector fields $X$, $Y$ and functions $f$ on $H^sG$.}

We also define a Riemannian metric $\langle \cdot, \cdot \rangle$ on $H^sG$ by letting
\begin{subequations}\label{metricdef}
\begin{align}\label{metricatid}
\langle (u, \rho), (v, \tau) \rangle_{(\text{id},0)}  \coloneqq &\; \langle u, v \rangle_{H^1} + \langle \rho, \tau \rangle_{L^2}
    \\ \nonumber
= &\; \int_{S^1} (uv + u_xv_x)dx + \int_{S^1} \rho \tau dx
\end{align}
at $(\id,0)$ and extending it to all of $H^sG$ by right-invariance:
\begin{equation}\label{metric}
\ska{X}{Y}_{(\varphi,f)}\coloneqq\ska{X(\varphi,f)\circ\varphi^{-1}}{Y(\varphi,f)\circ\varphi^{-1}}_{(\text{id},0)},
\end{equation}
\end{subequations}
where $X,Y$ are vector fields on $H^sG$. In the following, we will
write $\ska{\cdot}{\cdot}$ for $\ska{\cdot}{\cdot}_{(\text{id},0)}$.

It is a well-known fact that any Riemannian metric $\ska{\cdot}{\cdot}$ on a
finite-dimensional manifold $M$ induces a unique compatible \magenta{torsionless}
covariant derivative $\nabla$ on $M$ (\emph{the \red{L{e}vi}--Civita
connection}); $\nabla_XY$ is defined by
\begin{align}\label{levicivita}
2\ska{\nabla_XY}{Z}=&-\ska{[Y,X]}{Z}-\ska{X}{[Y,Z]}-\ska{Y}{[X,Z]}
\\ \nonumber
& +X\ska{Y}{Z}+Y\ska{Z}{X}-Z\ska{X}{Y},
\end{align}
where $X,Y,Z$ are vector fields on $M$. The bracket
$\ska{\cdot}{\cdot}$ establishes an isomorphism $T_mM \to T^*_mM$
for each $m\in M$, which guarantees the existence of
$\nabla_XY(m)$ for all $m$. However, this approach fails for
$H^sG$ since the metric defined by (\ref{metricdef}) is only a
\emph{weak} Riemannian metric (i.e. the topology induced by the
metric is weaker than the natural topology of any tangent space)
and therefore the metric does not establish an isomorphism between
the tangent space and its dual\magenta{, cf. \cite{CK02,EM70,L07}}.\\\indent
It is a first aim of this section
to establish that $\nabla$ as defined in \eqref{nabladef} defines a smooth
connection (i.e. $\Gamma$ defines a smooth spray) on $H^sG$ in the
sense of Banach manifolds (see \cite{Lang2}) and that $\langle
\cdot, \cdot \rangle$ is a compatible Riemannian metric.
\magenta{Note that
this connection is unique; this
can be deduced immediately from formula \eqref{levicivita}, as in the finite-dimensional case.}

In general, the Christoffel map for a Banach manifold is only defined locally.
Henceforth, we will implicitly use the natural smooth
identification
\begin{equation}\label{THsGidentification}
  TH^sG \simeq H^sG \times \left(H^{s}(S^1) \times H^{s -1}(S^1)\right)
\end{equation}
and view $\Gamma$ as a map from $H^sG$ to the space of bilinear
symmetric maps from $H^{s}(S^1) \times H^{s -1}(S^1)$ to itself.
Similarly, a vector field $X$ on $H^sG$ is viewed as a map $H^sG
\to H^{s}(S^1) \times H^{s -1}(S^1)$.
The identification (\ref{THsGidentification}) is given explicitly
as follows. The map $\varphi \mapsto (\varphi(0), \varphi(x) - x -
\varphi(0))$ is a diffeomorphism $\Diff(S^1) \to S^1 \times U^s$,
where
$$U^s := \{f \in H^s(S^1) | f(0) = 0, f_x > -1\}.$$
Since $U^s$ is an open subset of the closed linear subspace $E^s
:= \{f \in H^s(S^1) | f(0) = 0\} \subset H^s(S^1)$, this map
provides a local chart on $\Diff(S^1)$ with values in $I \times
U^s \subset \R \times E^s$ for any open subinterval $I \subset
S^1$. Moreover, using that $TS^1 \simeq S^1 \times \R$, we find
$$T\Diff(S^1) \simeq T(S^1 \times U^s) \simeq S^1 \times U^s \times \R \times E^s \simeq \Diff(S^1) \times H^s(S^1).$$
This yields the nontrivial part of (\ref{THsGidentification}).

For two Banach spaces $E, F$, we let
$\mathcal{L}^2_{\text{sym}}(E; F)$ denote the space of symmetric
bilinear maps from $E$ to $F$. For a manifold $M$,
$\mathcal{L}^2_{\text{sym}}(TM; F)$ denotes the bundle over $M$
with fiber $\mathcal{L}^2_{\text{sym}}(T_xM; F)$ over a point $x
\in M$.

\begin{prop}\label{2CHprop}
Let $s > 5/2$. Let $H^sG := H^{s}\Diff(\S)\circledS H^{s-1}(\S)$
and let $\Gamma$ be the Christoffel map defined in
(\ref{Gammadef}). Then $\Gamma$ defines a smooth spray on $H^sG$,
i.e., the map
\begin{align}\label{varphifGamma}
(\varphi,f)\mapsto \Gamma_{(\varphi,f)}\colon H^sG \to \mathcal{L}^2_{\text{\rm sym}}
\left(H^{s}(\S)\times H^{s-1}(\S);H^{s}(\S)\times H^{s-1}(\S)\right)
\end{align}
is smooth. Moreover, the metric $\ska{\cdot}{\cdot}$ defined by
(\ref{metricdef}) is a smooth (weak) Riemannian metric on $H^sG$,
i.e., the map
\begin{equation}\label{metricmap}
(\varphi,f) \mapsto \ska{\cdot}{\cdot}_{(\varphi,f)}: H^sG \to \mathcal L^2_{\text{\rm sym}}\left(T_{(\varphi,f)}H^sG;\R\right)
\end{equation}
is a smooth section of the bundle $\mathcal{L}^2_{\text{\rm
sym}}\left(TH^sG;\R\right)$.
Finally, the connection $\nabla$ and the metric $\ska{\cdot}{\cdot}$ are compatible in the sense that
\begin{equation}\label{compatible}
  X\ska{Y}{Z}=\ska{\nabla_XY}{Z} + \ska{Y}{\nabla_XZ}
\end{equation}
for all vector fields $X,Y, Z$ on $H^sG$.
\end{prop}
\proof
In order to establish smoothness of (\ref{varphifGamma}), it is sufficient to show that the following map is smooth:
\begin{align*}
& ((\varphi, f), w) \mapsto \Gamma_{(\varphi, f)}(w, w),
    \\
& H^sG \times [H^{s}(\S)\times H^{s-1}(\S)] \to H^{s}(\S)\times H^{s-1}(\S),
\end{align*}
where $w=(w_1,w_2)\in T_{(\varphi,f)}H^sG \simeq H^{s}(\S)\times H^{s-1}(\S)$ and
$$\Gamma_{(\varphi, f)}(w, w)  =
\begin{pmatrix}
\Gamma^0_{\text{id}}(w_1 \circ \varphi^{-1}, w_1 \circ \varphi^{-1}) -\frac{1}{2}A^{-1}\partial_x(w_2^2 \circ \varphi^{-1}) \\
  - (w_1\circ\varphi^{-1})_x w_2 \circ \varphi^{-1}
\end{pmatrix} \circ \varphi.
$$
We will show that the term
$-\frac{1}{2}\bigl(A^{-1}\partial_x(w_2^2 \circ
\varphi^{-1})\bigr) \circ \varphi$ makes a smooth contribution to
$\Gamma$; the other terms can be treated by similar arguments.

Consider the map
$$P:H^s\Diff(S^1) \times H^{s-1}(\S) \to H^s\Diff(S^1) \times H^{s}(\S)$$
defined by
\begin{align*}
P(\varphi, w) = \left(\varphi, \left(A^{-1}\partial_x(w^2 \circ
\varphi^{-1})\right)\circ \varphi\right).
\end{align*}
We write $P$ as the composition $P = \tilde{A}^{-1} \circ P_2 \circ P_1$, where the maps
\begin{align*}
&P_1:H^s\Diff(S^1) \times H^{s-1}(\S) \to H^s\Diff(S^1) \times H^{s-1}(\S),
    \\
&P_2: H^s\Diff(S^1) \times H^{s-1}(\S) \to H^s\Diff(S^1) \times H^{s-2}(\S),
    \\
& \tilde{A}:H^s\Diff(S^1) \times H^{s}(\S) \to H^s\Diff(S^1)
\times H^{s-2}(\S)
\end{align*}
are defined by
\begin{align*}
& P_1(\varphi, w)= (\varphi, w^2),
    \\
& P_2(\varphi, w) = \left(\varphi, (w \circ \varphi^{-1})_x \circ\varphi\right) = \left(\varphi, \frac{w_x}{\varphi_x}\right),
    \\
& \tilde{A}(\varphi, w) = (\varphi, \left(A( w \circ \varphi^{-1})\right) \circ\varphi)
= \left(\varphi, w - \frac{w_{xx}}{\varphi_x^2} + \frac{w_x\varphi_{xx}}{\varphi_x^3} \right).
\end{align*}
The maps $P_1, P_2$, and $\tilde{A}$ are smooth since $H^s(S^1)$ is a Banach algebra under pointwise multiplication for $s > 1/2$.
To show that $\tilde{A}^{-1}$ is smooth, we compute
$$D\tilde{A}(\varphi, w) = \begin{pmatrix} \id & 0 \\ * & \id - \frac{1}{\varphi_x^2}\partial_x^2 +
\frac{\varphi_{xx}}{\varphi_x^3}\partial_x \end{pmatrix}.$$ This
is, for each $(\varphi, w) \in H^s\Diff(S^1) \times H^{s}(\S)$, a
bijective bounded linear map $H^s(S^1) \times H^{s}(\S) \to
H^s(S^1) \times H^{s-2}(\S)$. The open mapping theorem implies
that its inverse is also bounded. The inverse mapping theorem now
implies that $\tilde{A}^{-1}$, and hence also $P$, is a smooth
map.

We next establish the smoothness of (\ref{metricmap}). It is
sufficient to show that the map
$$Q:H^sG \times \left[H^{s}(\S)\times H^{s-1}(\S)\right] \to \R,$$
defined by
\begin{equation}\label{localgdef}
  Q((\varphi,f), w) = \int_{\S}(w_1\circ\varphi^{-1})A(w_1\circ\varphi^{-1})dx+\int_{\S}(w_2\circ\varphi^{-1})^2 dx
\end{equation}
is smooth. The change of variables $y=\varphi^{-1}(x)$ yields
$$Q((\varphi,f), w)
= \int_{\S}\left(w_1^2\varphi_x + \frac{w_{1x}^2}{\varphi_x} + w_2^2\varphi_x\right)dy,$$
and written in this form the smoothness of $Q$ is clear.

It remains to verify (\ref{compatible}). Let $X_i, Y_i, Z_i$, $i =
1,2,$ denote the components of three vector fields $X,Y, Z$ on
$H^sG$. For $i = 1,2$, let $u_i =
X_i(\varphi,f)\circ\varphi^{-1}$,
$v_i=Y_i(\varphi,f)\circ\varphi^{-1}$,
$w_i=Z_i(\varphi,f)\circ\varphi^{-1}$. Let $\gamma(\epsilon)\in
H^sG$ be a curve such that $\gamma(0)= (\varphi,f)$ and
$\dot{\gamma}(0) = X(\varphi,f)$.

On the one hand,
\begin{align*}
(X\ska{Y}{Z})(\varphi,f) =\;& \frac{d}{d\epsilon}\bigg|_{\epsilon=0} \ska{Y(\gamma(\epsilon))}{Z(\gamma(\epsilon))}_{\gamma(\epsilon)}
    \\
=\;& \frac{d}{d\epsilon} \bigg|_{\epsilon=0}  \ska{Y_1(\gamma(\epsilon))\circ\gamma_1^{-1}}{Z_1(\gamma(\epsilon))\circ\gamma_1^{-1}}_{H^1}
    \\
& +\frac{d}{d \epsilon} \bigg|_{\epsilon=0} \ska{Y_2(\gamma(\epsilon))\circ\gamma_1^{-1}}{Z_2(\gamma(\epsilon))\circ\gamma_1^{-1}}_{L_2}.
\end{align*}

On the other hand,
\begin{align*}
\ska{\nabla_XY}{Z}_{(\varphi,f)}=&\ska{DY_1(\varphi,f)\cdot
X(\varphi,f)\circ\varphi^{-1}-\Gamma^0_\varphi(Y_1,X_1)\circ\varphi^{-1}}{w_1}_{H^1}
    \\
& +\frac{1}{2} \ska{(v_{2x}u_2+u_{2x}v_2)}{w_1}_{L_2}
    \\
&+\ska{DY_2(\varphi,f)\cdot
X(\varphi,f)\circ\varphi^{-1}+\frac{1}{2}(v_{1x}u_2+u_{1x}v_2)}{w_2}_{L_2}
\end{align*}
and
\begin{align*}
\ska{Y}{\nabla_XZ}_{(\varphi,f)}=&\ska{DZ_1(\varphi,f)\cdot
X(\varphi,f)\circ\varphi^{-1}-\Gamma^0_\varphi(Z_1,X_1)\circ\varphi^{-1}}{v_1}_{H^1}
    \\
&+\ska{\frac{1}{2}(w_{2x}u_2+u_{2x}w_2)}{v_1}_{L_2}
    \\
&+\ska{DZ_2(\varphi,f)\cdot
X(\varphi,f)\circ\varphi^{-1}+\frac{1}{2}(w_{1x}u_2+u_{1x}w_2)}{v_2}_{L_2}.
\end{align*}
The calculations in \cite{L07} for the CH equation show that
\bea
&&\left.\frac{d}{d\epsilon}\right|_{\epsilon =0} \ska{Y_1(\gamma(t))\circ\gamma_1^{-1}}{Z_1(\gamma(t))\circ\gamma_1^{-1}}_{H^1}=\nonumber\\
&&\qquad\ska{DY_1(\varphi,f)\cdot
X(\varphi,f)\circ\varphi^{-1}-\Gamma^0_\varphi(Y_1,X_1)\circ\varphi^{-1}}{w_1}_{H^1}\nonumber\\
&&\qquad+\ska{DZ_1(\varphi,f)\cdot
X(\varphi,f)\circ\varphi^{-1}-\Gamma^0_\varphi(Z_1,X_1)\circ\varphi^{-1}}{v_1}_{H^1},\nonumber
\eea
so it remains to check that
\bea
&&\left.\frac{d}{d\epsilon}\right|_{\epsilon =0}\ska{Y_2(\gamma(t))\circ\gamma_1^{-1}}{Z_2(\gamma(t))\circ\gamma_1^{-1}}_{L_2}\nonumber\\
&&\hspace{.5cm}=\ska{\frac{1}{2}(v_{2x}u_2+u_{2x}v_2)}{w_1}_{L_2}\nonumber\\
&&\hspace{.5cm}\quad+\ska{DY_2(\varphi,f)\cdot
X(\varphi,f)\circ\varphi^{-1}+\frac{1}{2}(v_{1x}u_2+u_{1x}v_2)}{w_2}_{L_2}\label{ddepsilonY2Z2}
\\
&&\hspace{.5cm}\quad+\ska{\frac{1}{2}(w_{2x}u_2+u_{2x}w_2)}{v_1}_{L_2}\nonumber\\
&&\hspace{.5cm}\quad+\ska{DZ_2(\varphi,f)\cdot
X(\varphi,f)\circ\varphi^{-1}+\frac{1}{2}(w_{1x}u_2+u_{1x}w_2)}{v_2}_{L_2}.\nonumber
\eea
Since
\begin{align*}
&\frac{d}{d\epsilon}\bigg|_{\epsilon =0} \ska{Y_2(\gamma(t))\circ\gamma_1^{-1}}{Z_2(\gamma(t))\circ\gamma_1^{-1}}_{L_2}
    \\
&\hspace{1cm} = \ska{DY_2(\varphi,f)\cdot X(\varphi,f)\circ\varphi^{-1}-v_{2x}u_1}{w_2}_{L_2}
    \\
&\hspace{1cm}\qquad+\ska{DZ_2(\varphi,f)\cdot X(\varphi,f)\circ\varphi^{-1}-w_{2x}u_1}{v_2}_{L_2},
\end{align*}
the condition in (\ref{ddepsilonY2Z2}) is equivalent to
\begin{align*}
\int_{\S}\bigg(& u_1v_{2x}w_2+u_1v_2w_{2x}+\frac{1}{2}u_2v_{2x}w_1+\frac{1}{2}u_{2x}v_2w_1+\frac{1}{2}u_2v_1w_{2x}
    \\
&+\frac{1}{2}u_{2x}v_1w_2+\frac{1}{2}u_2v_{1x}w_2+u_{1x}v_2w_2
+\frac{1}{2}u_2v_2w_{1x}\bigg)dx=0.
\end{align*}
Since the left-hand side is equal to
$$\int_{\S}\left(\frac{1}{2}\partial_x(u_2v_1w_2)+\frac{1}{2}\partial_x(u_2v_2w_1)+\partial_x(u_1v_2w_2)\right)dx=0,$$
we are done.
\endproof
\vspace{-.65cm}
\red{
\begin{rem}
%1. The Christoffel map $\Gamma$ is the infinite-dimensional analog
%of the Christoffel symbols $\Gamma^i_{jk}$ familiar from
%finite-dimensional differential geometry (see \cite{Lang2}; our
%$\Gamma$ is denoted by $B$ in \cite{Lang2}).
%
The crucial observation in the above proof is that $(P(w \circ
\varphi^{-1}))\circ\varphi$ is a rational expression in $w$,
$\varphi$, and their derivatives whenever $P$ is a differential
operator. This observation was already made on p. 154 of
\cite{EM70}.
\end{rem}}

A geodesic in $H^sG$ with respect to $\nabla$ is a $C^2$-curve
$(\varphi(t), f(t)) \in H^sG$ such that $\nabla_{(\varphi_t,f_t)}
(\varphi_t,f_t) = 0$, i.e.
\begin{equation}\label{geoeqn}
  (\varphi_{tt}, f_{tt}) = \Gamma_{(\varphi, f)}\bigl((\varphi_t,f_t), (\varphi_t,f_t)\bigr).
\end{equation}
Since the existence of a smooth connection on a Banach manifold
immediately yields the local existence and uniqueness of a
geodesic flow (see \cite{Lang2}), Proposition \ref{2CHprop}
implies the following result.

\begin{thm}\label{2CHgeodesicth}
Let $s > 5/2$. Then there exists an open interval $J$ centered at
$0$ and an open neighborhood $U$ of $(0,0) \in H^{s}(\S)\times
H^{s-1}(\S)$ such that for each $(u_0,\rho_0) \in U$ there exists
a unique solution $(\varphi, f) \in C^\infty\bigl(J, H^sG\bigr)$
of \eqref{geoeqn} with $(\varphi(0),f(0))=(\id,0)$ and
$(\varphi_t(0), f_t(0)) = (u_0,\rho_0)$. Furthermore, the solution
depends smoothly on the initial data in the sense that the local
flow $\Phi:J \times U \to H^sG$ defined by $\Phi(t, u_0, \rho_0) =
(\varphi(t; u_0, \rho_0), f(t; u_0, \rho_0)) $ is a smooth map.
\end{thm}

We write the Cauchy problem for 2CH in the form
\begin{align}\label{weak2CH}
  & \begin{pmatrix} u_t + uu_x \\ \rho_t + u\rho_x \end{pmatrix}
  = \begin{pmatrix}
  -A^{-1}\partial_x\biggl(u^2 + \frac{1}{2}u_x^2 + \frac{1}{2}\rho^2\biggr) \\ -\rho u_x
  \end{pmatrix},
    \\ \nonumber
  & (u(0), \rho(0)) = (u_0, \rho_0).
\end{align}
This formulation of 2CH is suitable for the formulation of weak
solutions. It follows from Theorem \ref{2CHgeodesicth} that the
2CH equation is locally well-posed in $H^{s}(S^1) \times
H^{s-1}(\S)$ for $s > 5/2$.

\begin{cor}[Local well-posedness in the $H^s$-category]\label{2CHlocalwellposedcor}
Suppose $s > 5/2$. Then for any $(u_0,\rho_0) \in H^{s}(S^1)
\times H^{s-1}(\S)$ there exists an open interval $J$ centered at
$0$ and a unique solution
\begin{equation}\label{urhosolution}
 (u, \rho) \in C\bigl(J, H^{s}(S^1) \times H^{s-1}(\S)\bigr) \cap C^1\bigl(J, H^{s-1}(S^1) \times H^{s-2}(\S)\bigr)
\end{equation}
of the Cauchy problem (\ref{weak2CH}) which depends continuously on the initial data $(u_0, \rho_0)$.
\end{cor}
\proof Theorem \ref{2CHgeodesicth} yields the existence of a
smooth curve $(\varphi(t), f(t)) \in H^sG$ such that
$(\varphi(0),f(0))=(\id,0)$ and $(\varphi_t(0), f_t(0)) =
(u_0,\rho_0)$. Define $(u(t), \rho(t))$ by equation
(\ref{Eulervel}). Then, $(u,\rho)$ has the regularity specified in
(\ref{urhosolution}) and depends continuously on $(u_0, \rho_0)$.
By right-invariance of $\Gamma$, the geodesic equation
(\ref{geoeqn}) can be written as
$$\begin{pmatrix} u_t + uu_x \\ \rho_t + u\rho_x \end{pmatrix}
= \Gamma_{(\id, 0)}((u,\rho), (u, \rho)).$$
This is equation (\ref{weak2CH}).
\endproof

\begin{rem}
The well-posedness result of Corollary \ref{2CHlocalwellposedcor}
can also be proved using Kato's semigroup approach (see
\cite{ELY07} for the case on the line).
\end{rem}

\subsection{The $C^n$-category}
The results of the previous subsection hold with the obvious
changes also in the $C^n$-category. Assuming $n \geq 2$, the
proofs are the same with $H^sG$ replaced with $C^nG$. In
particular, $\Gamma$ defines a smooth spray on $C^nG =
C^{n}\Diff(S^1) \circledS C^{n-1}(S^1)$ compatible with the metric
defined in (\ref{metricdef}). For the sake of brevity, we only
state the analog of Theorem \ref{2CHgeodesicth}.

\begin{thm}
Let $n \geq 2$. Then there exists an open interval $J$ centered at
$0$ and an open neighborhood $U$ of $(0,0) \in C^{n}(\S)\times
C^{n-1}(\S)$ such that for each $(u_0,\rho_0) \in U$ there exists
a unique solution $(\varphi, f) \in C^\infty\bigl(J, C^nG\bigr)$
of \eqref{geoeqn} with $(\varphi(0),f(0))=(\id,0)$ and
$(\varphi_t(0), f_t(0)) = (u_0,\rho_0)$. Furthermore, the solution
depends smoothly on the initial data in the sense that the local
flow $\Phi:J \times U \to C^nG$ defined by $\Phi(t, u_0, \rho_0) =
(\varphi(t; u_0, \rho_0), f(t; u_0, \rho_0))$ is a smooth map.
\end{thm}

\subsection{The smooth category}
We now want to extend the above results for 2CH to the space
$C^\infty G = C^\infty\Diff(S^1) \circledS C^\infty(S^1)$. Since
$C^\infty G$ is not a Banach manifold, the local existence and
uniqueness theorems for differential equations fail. We will
therefore take an indirect approach and first consider the local
geodesic flows on $H^sG$, $s > 5/2$. We will first show that the
domains of definition of these flows do not shrink to zero as $s
\to \infty$. By considering the limit as $s \to \infty$, the
existence of a smooth local geodesic flow on $C^\infty G$ will
then be established.

We will use the following blow-up result for 2CH.

\begin{prop}\label{blowupprop}
  Let $s > 5/2$. Let $(u_0, \rho_0) \in H^s(S^1) \times H^{s-1}(S^1)$ and let $T > 0$ be the maximal time of existence of the solution
  $$(u, \rho) \in C\bigl([0,T), H^s(S^1) \times H^{s-1}(\S)\bigr) \cap C^1\bigl([0,T), H^{s-1}(S^1) \times H^{s-2}(\S)\bigr)$$
of the Cauchy problem (\ref{weak2CH}). Then the solution $(u, \rho)$ blows up in finite time if and only if
\begin{equation}\label{blowupcriterion}
  \lim_{t \to T} \inf_{x \in S^1} \{u_x(t,x)\} = -\infty \quad \text{or} \quad \limsup_{t \to T} \{\|\rho_x(t)\|_{L^\infty}\} = \infty.
\end{equation}
\end{prop}
\proof
A proof for the equation obtained from (\ref{2CH}) by
replacing $\rho\rho_x$ with $-\rho\rho_x$ in the case on the line
is given in \cite{ELY07}; the same proof applies here.
\endproof

Let
$$\Phi_3: [0,T_3) \times U_3 \to H^3G,$$
where $T_3 > 0$ and $U_3 \subset H^{3}(\S)\times H^{2}(\S)$, be
the local geodesic flow on $H^3G$ whose existence is guaranteed by
Theorem \ref{2CHgeodesicth}. In the next proposition, we show that
the restriction of $\Phi_3$ to $H^{s}(\S)\times H^{s-1}(\S)$, $s
\geq 3$, defines a smooth flow on $H^sG$ for $t \in [0, T_3)$.
Thus, the flow on $H^sG$ exists for all $t \in [0, T_3)$ for any
$s \geq 3$.

\begin{prop}\label{2CHlocalflowprop}
Suppose $s > 3$ and let $\Phi_s$ denote the restriction of
$\Phi_3$ to $[0,T_3) \times U_s$, where $U_s = U_3 \cap
(H^{s}(\S)\times H^{s-1}(\S))$. Then $\Phi_s$ is a smooth local
flow of the geodesic equation \eqref{geoeqn} on $H^sG$, that is,
\begin{itemize}
\item[(a)] $\Phi_s$ is a smooth map from $[0,T_3) \times U_s$ to $H^sG$.
\item[(b)] For each $(u_0,\rho_0) \in U_s$, $\Phi_s(\cdot, u_0, \rho_0)$ is
a smooth solution of equation \eqref{geoeqn} on $[0, T_3)$ satisfying
$\Phi_s(0, u_0, \rho_0) = (\id,0)$ and $\partial_t\Phi_s(0, u_0, \rho_0) = (u_0,\rho_0)$.
\end{itemize}
\end{prop}
\proof Fix $(u_0,\rho_0) \in U_3$ and let $(u(t; u_0, \rho_0),
\rho(t; u_0, \rho_0))$ be the corresponding solution in
$H^{3}(S^1) \times H^{2}(\S)$ of the Cauchy problem
(\ref{weak2CH}). This solution is defined at least on $[0, T_3)$.
Since the criterion (\ref{blowupcriterion}) is independent of $s
\geq 3$, it follows from Proposition \ref{blowupprop} that if
$(u_0, \rho_0) \in U_s$ for some $s \geq 3$, then the curve $t
\mapsto (u(t;u_0, \rho_0), \rho(t;u_0,\rho_0))$ belongs to the
space
$$C\bigl([0,T_3), H^s(S^1) \times H^{s-1}(\S)\bigr) \cap C^1\bigl([0,T_3), H^{s-1}(S^1) \times H^{s-2}(\S)\bigr).$$
%Let $(\varphi, \rho) \in C^\infty ([0, T_3), H^3\Diff(S^1))$ be the geodesic $t \mapsto \Phi_3(t, u_0, \rho_0)$ defined on $[0, T_3)$.
Let $(\varphi,f)$ be the geodesic flow associated with the
solution $(u,\rho)$, defined on $[0,T_3)$.

Let $s > 3$. Suppose $(u_0, \rho_0) \in U_s$ and $\varphi \in
C^1([0, T_3), H^r\Diff(S^1))$ for some $r$ with $3 \leq r \leq
s-1$. We show that $\varphi \in C^1([0, T_3), H^{r+1}\Diff(S^1))$.
Using
$$\varphi_{tx}=(u_x\circ\varphi)\varphi_x, \qquad \varphi_{txx}=(u_{xx}\circ\varphi)\varphi_x^2+(u_x\circ\varphi)\varphi_{xx},$$
we find
$$\frac{d}{d t}\left(\frac{\varphi_{xx}}{\varphi_x}\right) =(u_{xx}\circ\varphi)\varphi_x.$$
Thus,
\begin{equation}\label{varphixx}
  \varphi_{xx}(t) = \varphi_x(t) \int_0^t (u_{xx}\circ\varphi)\varphi_x ds.
\end{equation}
Since $\varphi_x \in C^1([0, T_3), H^{r-1}(S^1))$ and $u_{xx} \in C\bigl([0,T_3), H^{s-2}(S^1)\bigr)$, equation (\ref{varphixx}) implies that
\begin{equation}\label{varphixxinC1}
  \varphi_{xx} \in C^1([0, T_3), H^{r-1}(S^1)).
\end{equation}
This implies that $\varphi \in C^1([0, T_3), H^{r+1}\Diff(S^1))$. Indeed,
\begin{align*}
\left\| \frac{\varphi(t) - \varphi(s)}{t-s} - u \circ \varphi\right\|_{H^{r+1}}^2
= & \left\| \frac{\varphi(t) - \varphi(s)}{t-s} - u \circ \varphi\right\|_{H^1}^2
    \\
& + \left\| \frac{\varphi_{xx}(t) - \varphi_{xx}(s)}{t-s} - (u \circ \varphi)_{xx} \right\|_{H^{r-1}}^2.
\end{align*}
As $t \to s$, the first term on the right-hand side vanishes
because $\varphi \in C^\infty([0, T_3), H^3\Diff(S^1))$ and the
second vanishes in view of (\ref{varphixxinC1}). Induction shows
that
\begin{equation}\label{varphiinC1}
  \varphi \in C^1([0, T_3), H^{s}\Diff(S^1)).
\end{equation}

We now show that in fact $(\varphi, f) \in C^\infty([0, T_3), H^{s}G)$. A computation shows that
\begin{equation}\label{ddtrhovarphi}
\frac{d}{d t}[(\rho\circ\varphi)\varphi_x]=[(\rho_{t}+u\rho_{x})\circ\varphi]\varphi_x+[(\rho u_x)\circ\varphi]\varphi_x=0.
\end{equation}
Thus, $f_t \varphi_x = (\rho\circ\varphi)\varphi_x = \rho_0$ and we infer that
\begin{equation}\label{ffromvarphi}
 f(t) = \rho_0 \int_0^t \frac{ds}{\varphi_x(s)}.
\end{equation}
It follows that
\begin{equation}\label{finC2}
f \in C^2([0, T_3), H^{s-1}(S^1)).
\end{equation}
Moreover, by Theorem \ref{2CHgeodesicth}, $(\varphi, f)$ is a
smooth solution of (\ref{geoeqn}) in $H^s\Diff(S^1) \times
H^{s-1}(S^1)$ for sufficiently small $t \geq 0$. Standard ODE
results show that the only way this solution can cease to exist
(Corollary IV.1.8 in \cite{Lang2}) is either that the condition
$\varphi_x > 0$ ceases to hold or that one of the norms
 \begin{equation}\label{norms}
  \left\|(\varphi_t, f_t)\right\|_{H^s(S^1) \times H^{s-1}(S^1)}, \quad \left\|\Gamma_{(\varphi, f)}
  \bigl((\varphi_t,f_t), (\varphi_t,f_t)\bigr) \right\|_{H^s (S^1) \times H^{s-1} (S^1)}
\end{equation}
blows up. But we know that $\varphi_x > 0$ on $[0, T_3)$ and
equations (\ref{varphiinC1}) and (\ref{finC2}) together with the
smoothness of $\Gamma$ imply that the norms in (\ref{norms})
remain bounded on $[0, T_3)$. This proves (b).

The standard ODE theorems on smooth dependence on initial data (Theorem IV.1.16 in \cite{Lang2}) imply (a).
\endproof

The Sobolev spaces $H^s(S^1)$ provide a Banach space approximation
of the Fr\'echet space $C^\infty(S^1)$ in the following sense.

\begin{defn}
  A {\it Banach space approximation} of a Fr\'echet space $X$ is a sequence of Banach spaces $\left(X_n, \|\cdot \|_n\right)_{n \geq 0}$ such that
  $$X_0 \supset X_1 \supset X_2 \supset \cdots \supset X \quad \text{and} \quad X = \cap_{n=0}^\infty X_n,$$
  where $\{\| \cdot \|\}_{n \geq 0}$ is a sequence of norms inducing the topology on $X$ such that
  $$\|x\|_0 \leq \|x\|_1 \leq \|x\|_2 \leq \cdots$$
  for all $x \in X$.
\end{defn}

The property of a Banach space approximation which is relevant for us is stated in the following lemma (a proof is given in \cite{EK09}).
\begin{lem}\label{banachapproxlemma}
  Let $X$ and $Y$ be Fr\'echet spaces with Banach space approximations $\{X_n\}_{n \geq 0}$ and $\{Y_n\}_{n \geq 0}$, respectively.
  Let $\Phi_0:U_0 \to V_0$ be a smooth map between two open subsets $U_0 \subset X_0$ and $V_0 \subset Y_0$. Let $U = U_0 \cap X$, $V = V_0 \cap Y$, and, for each $n \geq 0$,
  $$U_n = U_0 \cap X_n, \qquad V_n = V_0 \cap Y_n.$$
  Assume that, for each $n \geq 0$, the following properties are satisfied:
  \begin{itemize}
  \item[(1)] $\Phi_0(U_n) \subset V_n$,
  \item[(2)] the restriction $\Phi_0\big|_{U_n}:U_n \to V_n$ is a smooth map.
  \end{itemize}
  Then $\Phi_0(U) \subset V$ and the map $\Phi_0\big|_U:U \to V$ is smooth.
\end{lem}

Proposition \ref{2CHlocalflowprop} together with Lemma
\ref{banachapproxlemma} implies local well-posedness of the geodesic
flow on $C^\infty G$.

\begin{thm}\label{smooth2CHgeodesicth}
There exists an open interval $J$ centered at $0$ and an open
neighborhood $U$ of $(0,0) \in C^\infty(\S)\times C^\infty(\S)$
such that for each $(u_0,\rho_0) \in U$ there exists a unique
solution $(\varphi, f) \in C^\infty\bigl(J, C^\infty G\bigr)$ of
\eqref{geoeqn} satisfying $(\varphi(0),f(0))=(\id,0)$ and
$(\varphi_t(0), f_t(0)) = (u_0,\rho_0)$. Furthermore, the solution
depends smoothly on the initial data in the sense that the local
flow $\Phi:J \times U \to C^\infty G$ defined by $\Phi(t, u_0,
\rho_0) = (\varphi(t; u_0, \rho_0), f(t; u_0, \rho_0)) $ is a
smooth map.
\end{thm}

Since $C^\infty G$ is a Lie group with smooth multiplication and
$(u,\rho)=(\varphi_t\circ\varphi^{-1},f_t\circ\varphi^{-1})$,
\red{we have proved the first part of Theorem~\ref{thmmain}.}
%we immediately get the following result.

%\begin{cor}[Local well-posedness in the smooth category]
%There exists an open interval $J$ centered at $0$ and an open
%neighborhood $U$ of $(0,0) \in C^\infty(\S)\times C^\infty(\S)$
%such that for each $(u_0,\rho_0) \in U$ there exists a unique
%solution
%$$(u, \rho) \in C^\infty\bigl(J, C^\infty(S^1) \times C^\infty(S^1)\bigr)$$
%of (\ref{2CH}) with $(u(0), \rho(0)) = (u_0,\rho_0)$. Furthermore,
%the solution depends smoothly on the initial data in the sense
%that the local flow $\Phi:J \times U \to C^\infty(\S)\times
%C^\infty(\S)$ defined by $\Phi(t, u_0, \rho_0) = (u(t; u_0,
%\rho_0), \rho(t; u_0, \rho_0)) $ is a smooth map.
%\end{cor}

\section{The 2DP equation as a geodesic equation}\label{Sec2DP} \nequation
Most of the results for 2CH presented in the previous section have
direct counterparts in the case of 2DP; the main exception being
that the geodesic flow associated with 2DP is not induced by any
right-invariant metric. \red{(If this was the case, then,
choosing the second component to be equal to zero, we would obtain a
metric associated with the geodesic flow for DP which is not possible
as shown in \cite{ES10}.)}

\subsection{The $H^s$-category}
We define a bilinear operator $\Gamma_{(\text{id},0)}$
on\linebreak $H^{s}(\S)\times H^{s-1}(\S)$ by
\begin{subequations}\label{2DPGammadef}
\begin{equation}
\Gamma_{(\text{id},0)}((u,\rho), (v, \tau)) =
\begin{pmatrix}
\Gamma^0_{\text{id}}(u,v) - \frac{1}{2}A^{-1}(u_x \tau + v_x\rho) + A^{-1}\partial_x(\rho\tau) \\
  -(u_x\tau + v_x\rho)
\end{pmatrix},
\end{equation}
where $A = 1-\partial_x^2$ and
\begin{equation}
  \Gamma_{\text{id}}^0(u,v) = -\frac{3}{2} A^{-1}\partial_x\left(uv\right)
\end{equation}
\end{subequations}
is the Christoffel operator associated with the DP equation (cf. \cite{EK09}).
$\Gamma$ is extended to all of $H^sG$ by right-invariance, see Eq. (\ref{Gammarightinvariance}).
The corresponding covariant derivative $\nabla$ is defined by (\ref{nabladef}).
The proof of the following proposition is similar to that of Proposition \ref{2CHprop}.

\begin{prop}\label{2DPprop}
Let $s > 5/2$. Let $H^sG := H^{s}\Diff(\S)\circledS H^{s-1}(\S)$
and let $\Gamma$ be the 2DP Christoffel map defined in
(\ref{2DPGammadef}). Then $\Gamma$ defines a smooth spray on
$H^sG$, i.e., the map
\begin{align*}
(\varphi,f)\mapsto \Gamma_{(\varphi,f)}\colon H^sG \to \mathcal{L}^2_{\text{\rm sym}}(H^{s}(\S)\times H^{s-1}(\S);H^{s}(\S)\times
H^{s-1}(\S))
\end{align*}
is smooth.
\end{prop}

The existence of a smooth spray implies local existence and
uniqueness of the geodesic flow.

\begin{thm}\label{2DPgeodesicth}
Let $s > 5/2$. Let $\Gamma$ be the 2DP Christoffel map defined in
(\ref{2DPGammadef}). Then there exists an open interval $J$
centered at $0$ and an open neighborhood $U$ of $(0,0) \in
H^{s}(\S)\times H^{s-1}(\S)$ such that for each $(u_0,\rho_0) \in
U$ there exists a unique solution $(\varphi, f) \in
C^\infty\bigl(J, H^sG\bigr)$ of the geodesic equation
\eqref{geoeqn} satisfying $(\varphi(0),f(0))=(\id,0)$ and
$(\varphi_t(0), f_t(0)) = (u_0,\rho_0)$. Furthermore, the solution
depends smoothly on the initial data in the sense that the local
flow $\Phi:J \times U \to H^sG$ defined by $\Phi(t, u_0, \rho_0) =
(\varphi(t; u_0, \rho_0), f(t; u_0, \rho_0)) $ is a smooth map.
\end{thm}

We write the Cauchy problem for 2DP in the form
\begin{align}\label{weak2DP}
  & \begin{pmatrix} u_t + uu_x \\ \rho_t + u\rho_x \end{pmatrix}
  = \begin{pmatrix}
    -A^{-1}\left(\left(\frac{3}{2}u^2 - \rho^2\right)_x + \rho u_x \right) \\
  -2\rho u_x
  \end{pmatrix},
    \\ \nonumber
  & (u(0), \rho(0)) = (u_0, \rho_0).
\end{align}
It follows from Theorem \ref{2DPgeodesicth} that 2DP is locally
well-posed in $H^{s}(S^1) \times H^{s-1}(\S)$ for $s > 5/2$.

\begin{cor}[Local well-posedness in the $H^s$-category]
Suppose $s > 5/2$. Then for any $(u_0,\rho_0) \in H^{s}(S^1)
\times H^{s-1}(\S)$ there exists an open interval $J$ centered at
$0$ and a unique solution
\begin{equation*}
 (u, \rho) \in C\bigl(J, H^{s}(S^1) \times H^{s-1}(\S)\bigr) \cap C^1\bigl(J, H^{s-1}(S^1) \times H^{s-2}(\S)\bigr)
\end{equation*}
of the Cauchy problem (\ref{weak2DP}) which depends continuously on the initial data $(u_0, \rho_0)$.
\end{cor}
\proof\magenta{Let $(\varphi(t), f(t)) \in H^sG$ be the smooth curve with
$(\varphi(0),f(0))=(\id,0)$ and $(\varphi_t(0), f_t(0)) = (u_0,\rho_0)$ obtained in Theorem \ref{2DPgeodesicth} and define $(u(t), \rho(t)):=(\varphi_t(t),f_t(t))\circ\varphi^{-1}(t)$. Then, $(u,\rho)$ has the regularity specified in
the corollary and depends continuously on $(u_0, \rho_0)$.
By right-invariance of the 2DP Christoffel map $\Gamma$, the geodesic equation
$(\varphi_{tt},f_{tt})=\Gamma_{(\varphi,f)}((\varphi_t,f_t),(\varphi_t,f_t))$ can be written as
$$\begin{pmatrix} u_t + uu_x \\ \rho_t + u\rho_x \end{pmatrix}
= \Gamma_{(\id, 0)}((u,\rho), (u, \rho)).$$
This is equation (\ref{weak2DP}).}
\endproof

\subsection{The $C^n$-category}
The results of the previous subsection hold with the obvious changes also in the $C^n$-category, $n \geq 2$.

\subsection{The smooth category}
We have the following blow-up result for 2DP; the proof is similar
to that of Proposition \ref{blowupprop}.

\begin{prop}\label{2DPblowupprop}
  Let $s > 5/2$. Let $(u_0, \rho_0) \in H^s(S^1) \times H^{s-1}(S^1)$ and let $T > 0$ be the maximal time of existence of the solution
  $$(u, \rho) \in C\bigl([0,T), H^s(S^1) \times H^{s-1}(\S)\bigr) \cap C^1\bigl([0,T), H^{s-1}(S^1) \times H^{s-2}(\S)\bigr)$$
of the Cauchy problem (\ref{weak2DP}). Then the solution $(u, \rho)$ blows up in finite time if and only if
\begin{equation*}
  \lim_{t \to T} \inf_{x \in S^1} \{u_x(t,x)\} = -\infty \quad \text{or} \quad \limsup_{t \to T} \{\|\rho_x(t)\|_{L^\infty}\} = \infty.
\end{equation*}
\end{prop}

Let
$$\Phi_3: [0,T_3) \times U_3 \to H^3G,$$
where $T_3 > 0$ and $U_3 \subset H^{3}(\S)\times H^{2}(\S)$, be
the local geodesic flow on $H^3G$ whose existence is guaranteed by
Theorem \ref{2DPgeodesicth}.

\begin{prop}\label{2DPlocalflowprop}
Suppose $s > 3$ and let $\Phi_s$ denote the restriction of
$\Phi_3$ to $[0,T_3) \times U_s$, where $U_s = U_3 \cap
(H^{s}(\S)\times H^{s-1}(\S))$. Let $\Gamma$ be the 2DP
Christoffel map defined in (\ref{2DPGammadef}). Then $\Phi_s$ is a
smooth local flow of the geodesic equation \eqref{geoeqn} on
$H^sG$, that is,
\begin{itemize}
\item[(a)] $\Phi_s$ is a smooth map from $[0,T_3) \times U_s$ to
$H^sG$. \item[(b)] For each $(u_0,\rho_0) \in U_s$, $\Phi_s(\cdot,
u_0, \rho_0)$ is a smooth solution of equation \eqref{geoeqn} on
$[0, T_3)$ satisfying $\Phi_s(0, u_0, \rho_0) = (\id,0)$ and
$\partial_t\Phi_s(0, u_0, \rho_0) = (u_0,\rho_0)$.
\end{itemize}
\end{prop}
\proof The proof is identical to that of Proposition
\ref{2CHlocalflowprop} except that equation (\ref{ffromvarphi})
must be replaced with
\begin{equation}\label{2DPffromvarphi}
  f(t) = \rho_0 \int_0^t \frac{ds}{\varphi_x^2(s)}.
\end{equation}
Equation (\ref{2DPffromvarphi}) is proved by noting that
$$\frac{d}{d t}\left[(\rho\circ\varphi)\varphi_x^2\right]=[(\rho_{t}+u\rho_{x})\circ\varphi]\varphi_x^2+2[(\rho u_x)\circ\varphi]\varphi_x^2=0,$$
and so $f_t \varphi_x^2 = (\rho\circ\varphi)\varphi_x^2 = \rho_0$.
\endproof

We find the following well-posedness results.

\begin{thm}
Let $\Gamma$ be the 2DP Christoffel map. There exists an open
interval $J$ centered at $0$ and an open neighborhood $U$ of
$(0,0) \in C^\infty(\S)\times C^\infty(\S)$ such that for each
$(u_0,\rho_0) \in U$ there exists a unique solution $(\varphi, f)
\in C^\infty\bigl(J, C^\infty G\bigr)$ of the geodesic equation
\eqref{geoeqn} satisfying $(\varphi(0),f(0))=(\id,0)$ and
$(\varphi_t(0), f_t(0)) = (u_0,\rho_0)$. Furthermore, the solution
depends smoothly on the initial data in the sense that the local
flow $$\Phi:J \times U \to C^\infty G,\quad\Phi(t, u_0, \rho_0) =
(\varphi(t; u_0, \rho_0), f(t; u_0, \rho_0))$$ is a smooth map.
\end{thm}

\red{By the same arguments as in the previous section, this proves the
second part of Theorem~\ref{thmmain}. Hence the proof of Theorem~\ref{thmmain} is completed.}
%\begin{cor}[Local well-posedness in the smooth category]
%There exists an open interval $J$ centered at $0$ and an open
%neighborhood $U$ of $(0,0) \in C^\infty(\S)\times C^\infty(\S)$
%such that for each $(u_0,\rho_0) \in U$ there exists a unique
%solution
%$$(u, \rho) \in C^\infty\bigl(J, C^\infty(S^1) \times C^\infty(S^1)\bigr)$$
%of (\ref{2DP}) with $(u(0), \rho(0)) = (u_0,\rho_0)$. Furthermore,
%the solution depends smoothly on the initial data in the sense
%that the local flow $\Phi:J \times U \to C^\infty(\S)\times
%C^\infty(\S)$ given by $\Phi(t, u_0, \rho_0) = (u(t; u_0, \rho_0),
%\rho(t; u_0, \rho_0)) $ is a smooth map.
%\end{cor}

\section{The sectional curvature for the 2CH equation}\label{SecCurvature} \nequation
We have showed that both 2CH and 2DP are geodesic equations on
$H^sG = H^s\Diff(S^1) \circledS H^{s-1}(S^1)$ with respect to a
smooth affine connection. The existence of a smooth connection
$\nabla$ on a Banach manifold immediately implies the existence of
a smooth curvature tensor $R$ defined by
$$R(X,Y)Z=\nabla_X\nabla_YZ-\nabla_Y\nabla_XZ-\nabla_{[X,Y]}Z,$$
where $X,Y, Z$ are vector fields on $H^sG$ (cf. \cite{Lang2}). In
the case of 2CH, since there exists a metric $\langle \cdot, \cdot
\rangle$, we can also define an (unnormalized) sectional curvature
tensor $S$ by\footnote{Recall that the sectional curvature
$\text{Sec}(\sigma)$ of a subspace $\sigma$ spanned by two tangent
vectors $u$ and $v$ is defined by
$$\text{Sec}(\sigma) = \frac{\langle R(u,v)v, u\rangle}{\langle u, u \rangle \langle v, v\rangle - \langle u, v \rangle^2}.$$}
$$S(X,Y) := \langle R(X,Y)Y, X \rangle.$$
In this section, we will derive a convenient formula for $S$ and
use it to determine large subspaces of positive curvature for the
2CH equation.

We will work in the $H^s$-category; similar results are valid with
$H^sG$ replaced with $C^nG$. In view of the right-invariance of
$\nabla$, it is enough to consider the curvature at the identity
$(\id, 0)$. We will write $\Gamma$ for $\Gamma_{(\id, 0)}$.

\begin{prop} Let $s > 5/2$. Let $R$ be the curvature tensor on
$H^sG$ associated with the 2CH equation. Then $S(u,v) :=\ska{R(u,v)v}{u}$
is given at the identity by
\begin{equation}\label{Suvexpression}
  S(u,v)=\ska{\Gamma(u,v)}{\Gamma(u,v)}-\ska{\Gamma(u,u)}{\Gamma(v,v)}, \qquad u,v\in T_{(\id,0)}H^sG.
\end{equation}
\end{prop}
\proof
Let $U, V, W \in T_{p}H^sG$ be three tangent vectors at a point $p \in H^sG$.
The curvature tensor $R$ is given locally by \cite{Lang2}
\begin{align*}
R_p(U, V)W = & D_1\Gamma_{p}(W, U)V - D_1\Gamma_{p}(W, V)U
         \\
 &  + \Gamma_{p}(\Gamma_{p}(W, V), U) - \Gamma_{p}(\Gamma_{p}(W, U), V)
\end{align*}
where $\Gamma$ is the 2CH Christoffel map defined in
(\ref{Gammadef}) and $D_1$ denotes differentiation with respect to
$p$:
$$D_1\Gamma_{p}(W, U)V = \frac{d}{d\epsilon}\bigg|_{\epsilon = 0} \Gamma_{p + \epsilon V}(W, U).$$
Let $\mathfrak{g}_s := T_{(\id,0)}H^sG$. Let $u = (u_1, u_2)$, $v
= (v_1, v_2)$, and $w = (w_1, w_2)$ be three vectors in
$\mathfrak{g}_s \simeq H^s(S^1) \times H^{s-1}(S^1)$. Using the
identity
$$\left.\frac{d}{d\eps}\right|_{\eps=0} u_1\circ(\id+\eps v_1)^{-1} =-u_{1x}v_1$$
a long but straightforward computation shows that
$$D_1\Gamma(w,u)v=-\Gamma(w_xv_1,u)-\Gamma(u_xv_1,w)+\Gamma(w,u)_xv_1.$$
Thus,
\begin{align}\nonumber
S(u,v) =&\; \ska{\Gamma(\Gamma(v,v),u)}{u}-\ska{\Gamma(\Gamma(v,u),v)}{u}
    \\ \label{Suv}
& + \langle \Gamma(v,u)_xv_1, u \rangle - \langle \Gamma(v,v)_xu_1, u \rangle
    \\ \nonumber
&+\ska{-\Gamma(v_xv_1,u)-\Gamma(v,u_xv_1)+2\Gamma(v_xu_1,v)}{u}.
\end{align}

We define a bilinear operator $B = (B_1, B_2):\mathfrak{g}_s \times \mathfrak{g}_s \to \mathfrak{g}_s$ by
$$\begin{pmatrix}
  B_1(u,v) \\
  B_2(u,v)
\end{pmatrix}
= \begin{pmatrix}
  -A^{-1}(2v_{1x}Au_1+v_1Au_{1x}+u_2v_{2x}) \\
  -(u_2v_1)_x
\end{pmatrix}.
$$
Then $B$ satisfies $\ska{B(u,v)}{w}=\ska{u}{[v,w]}$ where $[v,w] = w_xv - v_xw$ and
\begin{equation}\label{Gamma2CH}
\Gamma(u,v)=\frac{1}{2}\left[\begin{pmatrix}
  (u_1v_1)_x \\
  u_{2x}v_1+v_{2x}u_1
\end{pmatrix}
+B(u,v)+B(v,u)\right].
\end{equation}
Let $\Gamma_1$ and $\Gamma_2$ denote the two components of $\Gamma$.
With this notation, the first four terms on the right-hand side of (\ref{Suv}) equal
\begin{align*}
&\frac{1}{2}\ska{\begin{pmatrix}
  (\Gamma_1(v,v)u_1)_x \\
  \Gamma_2(v,v)_xu_1+u_{2x}\Gamma_1(v,v)
  \end{pmatrix}
  +B(\Gamma(v,v),u)+B(u,\Gamma(v,v))}{u}
    \\
&-\frac{1}{2}\ska{\begin{pmatrix}
  (\Gamma_1(v,u)v_1)_x \\
  \Gamma_2(v,u)_xv_1+v_{2x}\Gamma_1(v,u)
  \end{pmatrix}
  +B(\Gamma(v,u),v)+B(v,\Gamma(v,u))}{u}
    \\
&+ \langle \Gamma(v,u)_xv_1, u \rangle - \langle \Gamma(v,v)_xu_1, u \rangle.
\end{align*}
We rewrite this expression as
$$\frac{1}{2}\ska{[v,\Gamma(v,u)]}{u}+\ska{u}{[\Gamma(v,v),u]}
-\frac{1}{2}\ska{\Gamma(v,u)}{[v,u]}-\frac{1}{2}\ska{v}{[\Gamma(v,u),u]},$$
which in turn equals
\begin{align*}
& \ska{\Gamma(u,v)}{\Gamma(u,v)}-\ska{\Gamma(u,u)}{\Gamma(v,v)}
    \\
&+\ska{\begin{pmatrix}
  u_{1x}u_1 \\
  u_{2x}u_1
\end{pmatrix}}{\Gamma(v,v)} -\ska{\begin{pmatrix}
  u_{1x}v_1 \\
  u_{2x}v_1
\end{pmatrix}}{\Gamma(u,v)}.
\end{align*}
Hence, equation (\ref{Suv}) becomes
\begin{align} \nonumber
S(u,v)=&\;\ska{\Gamma(u,v)}{\Gamma(u,v)}-\ska{\Gamma(u,u)}{\Gamma(v,v)}
    \\ \label{Suv2}
&-\ska{\begin{pmatrix}
  u_{1x}v_1 \\
  u_{2x}v_1
\end{pmatrix}}{\Gamma(u,v)}+\ska{
\begin{pmatrix}
  u_{1x}u_1 \\
  u_{2x}u_1
\end{pmatrix}}{\Gamma(v,v)}
    \\ \nonumber
&+\ska{-\Gamma(v_xv_1,u)-\Gamma(v,u_xv_1)+2\Gamma(v_xu_1,v)}{u}
\end{align}
We claim that the sum of the last three terms on the
right-hand side of (\ref{Suv2}) is zero. Indeed, using the
expression
\begin{equation}\label{Gammauv}
\Gamma(u,v)=\left(%
\begin{array}{c}
  \Gamma^0(u_1,v_1)-\frac{1}{2}A^{-1}(u_2v_2)_x \\
  -\frac{1}{2}(u_{1x}v_2+v_{1x}u_2) \\
\end{array}%
\right)
\end{equation}
for $\Gamma$, integration by parts shows that the terms in
(\ref{Suv2}) involving $\Gamma^0$ cancel. A somewhat tedious
computation involving further integration by parts shows that the
remaining terms also vanish. This proves (\ref{Suvexpression}).
\endproof

A formula analogous to (\ref{Suvexpression}) for the CH equation
was derived in \cite{LMP}: If $S_{CH}(u_1,v_1)$ denotes the
unnormalized sectional curvature on $H^s\Diff(S^1)$ associated
with the CH equation, then
$$S_{CH}(u_1,v_1) = \ska{\Gamma^0(u_1,v_1)}{\Gamma^0(u_1,v_1)}-\ska{\Gamma^0(u_1,u_1)}{\Gamma^0(v_1,v_1)},$$
for all $u_1, v_1 \in T_\id H^s\Diff(S^1)$.
It was also shown in \cite{LMP} that
\begin{align}\label{SCHcoscos}
S_{CH}&(\cos kx,\cos lx)
    \\ \nonumber
&= \frac{1}{8}\left(\frac{(1+\frac{1}{2}kl)^2}{1+(k-l)^2}(k-l)^2+
\frac{(1-\frac{1}{2}kl)^2}{1+(k+l)^2}(k+l)^2\right)>0,
\end{align}
whenever $k, l \in \{ 2\pi, 4\pi, \dots \}$, $k \neq l$,
establishing the existence of a large subspace of positive
curvature for CH. Since
\begin{equation}\label{SSCHrelation}
S\biggl(\begin{pmatrix}
  u_1 \\
  0
  \end{pmatrix},
\begin{pmatrix}
  v_1 \\
  0
  \end{pmatrix}\biggr) = S_{CH}(u_1,v_1),
\end{equation}
we conclude that the same example yields an infinite-dimensional
subspace of positive curvature for 2CH. In the next proposition,
we investigate the curvature of $\Diff(S^1) \circledS
\mathcal{F}(S^1)$ in directions which are nontrivial along the
second component.

\begin{prop}\label{curvprop}
 Let $s > 5/2$. Let $S(u,v) :=\ska{R(u,v)v}{u}$ be the
 unnormalized sectional curvature on $H^sG$ associated with the 2CH equation.
Then
$$S(u,v) > 0$$
for all vectors $u, v \in T_{(\id, 0)}H^sG$, $u \neq v$, of the form
\begin{equation}\label{uvcoscos}
u = \begin{pmatrix} \cos k_1 x \\ \cos k_2x \end{pmatrix}, \quad
v = \begin{pmatrix} \cos l_1 x \\ \cos l_2x \end{pmatrix}, \qquad k_1, k_2, l_1, l_2 \in \{2\pi, 4\pi, \dots\}.
\end{equation}
Moreover, the sectional curvature $\text{\upshape Sec}(u,v)$ satisfies
\begin{equation}\label{Secuv}
\text{\upshape Sec}(u,v) :=\frac{S(u,v)}{\langle u, u \rangle
\langle v, v\rangle - \langle u, v \rangle^2}  \geq \frac{1}{8}
\end{equation}
for all vectors $u, v \in T_{(\id, 0)}H^sG$, $u \neq v$, of the form
\begin{equation}\label{uvfirstzero}
u = \begin{pmatrix} 0 \\ \cos k_2x \end{pmatrix}, \quad
v = \begin{pmatrix} 0 \\ \cos l_2x \end{pmatrix}, \qquad k_2, l_2 \in \{2\pi, 4\pi, \dots\}.
\end{equation}
\end{prop}
\proof
In view of (\ref{Suvexpression}), we have
\begin{align*}
S(u,v) =&\int_{\S}\Gamma_1(u,v)A\Gamma_1(u,v)dx+\int_{\S}\Gamma_2(u,v)^2dx
    \\
&-\int_{\S}\Gamma_1(u,u)A\Gamma_1(v,v)dx-\int_{\S}\Gamma_2(u,u)\Gamma_2(v,v)dx.
\end{align*}
Using the expression (\ref{Gammauv}) for $\Gamma(u,v)$ and integrating by parts, we find
\begin{align}\label{SuvsumIj}
S(u,v) = S_{CH}(u_1,v_1) + \sum_{j=1}^4 I_j,
\end{align}
where
\begin{align*}
& I_1 = \frac{1}{4}\int_{\S}(u_2v_2)_xA^{-1}(u_2v_2)_xdx
    \\
& I_2 = -\frac{1}{4}\int_{\S}(u_2^2)_xA^{-1}(v_2^2)_xdx,
    \\
& I_3 = \frac{1}{2}\int_{\S} \left[\Gamma^0(u_1,u_1)(v_2^2)_x + \Gamma^0(v_1,v_1)(u_2^2)_x - 2\Gamma^0(u_1,v_1)(u_2v_2)_x\right] dx,
    \\
&I_4 = \frac{1}{4}\int_{\S}(u_{1x}^2v_2^2 + v_{1x}^2u_2^2)dx - \frac{1}{2} \int_{\S} u_{1x}u_2v_{1x}v_2dx.
\end{align*}
Now suppose $u$ and $v$ have the form specified in (\ref{uvcoscos}).
Then the terms $\{I_j\}_1^4$ can be computed explicitly using the trigonometric identities
\begin{align*}
&\cos\alpha\cos\beta=\frac{1}{2}(\cos(\alpha-\beta)+\cos(\alpha+\beta)),
    \\
&\sin\alpha\sin\beta=\frac{1}{2}(\cos(\alpha-\beta)-\cos(\alpha+\beta)),
    \\
&\sin\alpha\cos\beta=\frac{1}{2}(\sin(\alpha-\beta)+\sin(\alpha+\beta)),
\end{align*}
the relations
\begin{align*}
 A^{-1}\cos\alpha x= &\;\frac{1}{1+\alpha^2}\cos \alpha x, \qquad \alpha \in \R,
    \\
 \int_0^1\cos(\alpha x)\cos(\beta x) dx= &\;\frac{1}{2}\left(\delta_{\alpha,\beta}+\delta_{\alpha,-\beta}\right), \qquad \alpha,\beta\in 2\pi\Z,
    \\
 \int_0^1\sin(\alpha x)\sin(\beta x)dx = &\; \frac{1}{2}\left(\delta_{\alpha,\beta} - \delta_{\alpha,-\beta}\right), \qquad \alpha,\beta\in 2\pi\Z,
    \\
 \int_0^1\cos(\alpha x)\sin(\beta x)dx = &\;0, \qquad \alpha,\beta\in 2\pi\Z,
\end{align*}
and the identity
\begin{align*}
\Gamma^0&(\cos \alpha x,\cos \beta x)
    \\
&=\partial_x \bigg[-\frac{\frac{1}{2}(1-\frac{1}{2}\alpha \beta)}{1+(\alpha+ \beta)^2}\cos(\alpha+ \beta)x
-\frac{\frac{1}{2}(1+\frac{1}{2}\alpha \beta)}{1+(\alpha-\beta)^2}\cos(\alpha-\beta)x\bigg],
    \\
&\hspace{9cm} \alpha, \beta \in 2\pi \Z.
\end{align*}
We find
\begin{align} \nonumber
 I_1 = &\; \frac{1}{32}\left(\frac{(k_2-l_2)^2}{1+(k_2-l_2)^2}+\frac{(k_2+l_2)^2}{1+(k_2+l_2)^2}\right),
    \\ \nonumber
 I_2 = &-\frac{1}{8}\frac{k_2^2}{1+(2k_2)^2}\delta_{k_2,l_2},
    \\ \nonumber
 I_3 = &\;\frac{1}{8}\frac{(1-\frac{1}{2}k_1l_1)(k_1+l_1)^2}{1+(k_1+l_1)^2}
\left(\delta_{k_1+l_1,k_2-l_2}+\delta_{k_1+l_1,l_2-k_2}+\delta_{k_1+l_1,k_2+l_2}\right)
    \\ \nonumber
&+\frac{1}{8}\frac{(1+\frac{1}{2}k_1l_1)(k_1-l_1)^2}{1+(k_1-l_1)^2}
    \\ \label{I1234}
& \quad \times (\delta_{k_1-l_1,k_2-l_2}+\delta_{k_1-l_1,l_2-k_2} +\delta_{k_1-l_1,k_2+l_2}+\delta_{l_1-k_1,k_2+l_2})
    \\ \nonumber
&- \frac{k_1^2}{4}\frac{1-\frac{1}{2}k_1^2}{1+(2k_1)^2}\delta_{k_1,l_2}
- \frac{l_1^2}{4}\frac{1-\frac{1}{2}l_1^2}{1+(2l_1)^2}\delta_{k_2,l_1},
    \\ \nonumber
I_4 = &\; \frac{1}{16}k_1^2\left(1-\frac{1}{2}\delta_{k_1,l_2}\right)+\frac{1}{16}l_1^2\left(1-\frac{1}{2}\delta_{l_1,k_2}\right)
    \\ \nonumber
&-\frac{1}{16}k_1l_1\bigl(\delta_{k_1-l_1,k_2-l_2}+\delta_{k_1-l_1,l_2-k_2}+\delta_{k_1-l_1,k_2+l_2}
+\delta_{l_1-k_1,k_2+l_2}
    \\ \nonumber
& \qquad \qquad \;\; -\delta_{k_1+l_1,k_2-l_2}-\delta_{k_1+l_1,l_2-k_2}-\delta_{k_1+l_1,k_2+l_2}\bigr).
\end{align}
Together with expression (\ref{SCHcoscos}) for $S_{CH}(u_1,v_1)$
this yields an expression for $S(u,v)$ in terms of $k_1, k_2, l_1,
l_2$. The sum of the negative terms in this expression can be
estimated as follows:
\begin{align}\nonumber
& -\frac{1}{8}\frac{k_2^2}{1+(2k_2)^2}\delta_{k_2,l_2}
    \\ \label{negativeterms}
& -\frac{1}{16}k_1l_1\frac{(k_1+l_1)^2}{1+(k_1+l_1)^2}
\left(\delta_{k_1+l_1,k_2-l_2}+\delta_{k_1+l_1,l_2-k_2}+\delta_{k_1+l_1,k_2+l_2}\right)
    \\ \nonumber
&-\frac{1}{16}k_1l_1\bigl(\delta_{k_1-l_1,k_2-l_2}+\delta_{k_1-l_1,l_2-k_2}+\delta_{k_1-l_1,k_2+l_2}
+\delta_{l_1-k_1,k_2+l_2}\bigr)
    \\ \nonumber
\geq & -\frac{1}{32} - \frac{k_1l_1}{16} - \frac{k_1l_1}{16},
\end{align}
because at most one delta function within each bracket can give a
nonzero contribution for a given set of values of $k_1, k_2, l_1,
l_2 \in 2\pi\N$.

On the other hand, the term $S_{CH}(u_1,v_1)$ contributes to $S(u,v)$ the positive term
\begin{equation}\label{positiveterm}
  \frac{1}{8}\frac{(1-\frac{1}{2}k_1l_l)^2}{1+(k_1+l_1)^2}(k_1+l_1)^2,
\end{equation}
and the sum of the right-hand side of (\ref{negativeterms}) and (\ref{positiveterm}) is positive:
\begin{align*}
\frac{1}{8}\frac{(1-\frac{1}{2}k_1l_l)^2}{1+(k_1+l_1)^2}&(k_1+l_1)^2
-\frac{1}{32} - \frac{k_1l_1}{8}
    \\
&\geq \frac{1}{16}\left(1-\frac{1}{2}k_1l_l\right)^2 -\frac{1}{32}
- \frac{k_1l_1}{8}
    \\
&= \frac{k_1^2l_1^2}{16}\left[\frac{1}{k_1^2l_1^2} -
\frac{1}{k_1l_l} + \frac{1}{4} - \frac{1}{2 k_1^2 l_1^2} -
\frac{2}{k_1l_2}\right] > 0,
\end{align*}%
where we used that $k_1, l_1 \geq 2\pi$. This shows that $S(u,v) > 0$.

In remains to prove (\ref{Secuv}). Suppose $u_1 = v_1 = 0$ and $u_2 \neq v_2$. It follows from (\ref{SuvsumIj}) and (\ref{I1234}) that
\begin{align} \nonumber
S\biggl(\begin{pmatrix} 0 \\ u_2\end{pmatrix}&, \begin{pmatrix} 0 \\ v_2 \end{pmatrix}\biggr)
=  I_1 + I_2
    \\ \label{S0u20v2}
= &\; \frac{1}{32}\left(\frac{(k_2-l_2)^2}{1+(k_2-l_2)^2}+\frac{(k_2+l_2)^2}{1+(k_2+l_2)^2}\right)
-\frac{1}{8}\frac{k_2^2}{1+(2k_2)^2}\delta_{k_2,l_2}
    \\ \nonumber
\geq & \;\frac{1}{64} + \frac{1}{64},
\end{align}
where we used that $k_2 \neq l_2$.
On the other hand, for this choice of $u$ and $v$,
$$\langle u, v\rangle = \frac{1}{2} \delta_{k_2, l_2},$$
and hence
\begin{align}\label{uvareaelement}
\langle u, u \rangle \langle v, v\rangle - \langle u, v \rangle^2
= \frac{1}{4}.
\end{align}
Equations (\ref{S0u20v2}) and (\ref{uvareaelement}) yield (\ref{Secuv}).
\endproof

\begin{rem}
Although Proposition \ref{curvprop} establishes the existence of a
large subspace of positive curvature, there are also directions
for 2CH of strictly negative curvature. Indeed, it is shown in
\cite{LMP} that there exist directions of strictly negative
sectional curvature for the CH equation. In view of
(\ref{SSCHrelation}), this implies that 2CH also admits directions
of negative curvature.
\end{rem}

\appendix
\section{Comparison with the rotating rigid body}\label{SecRigidBody}\nequation
\renewcommand{\theequation}{A.\arabic{equation}}
In this appendix, the geometric interpretations of 2CH, CH, and
the rotating rigid body are compared in an attempt to emphasize
some unifying features of the approach pioneered by Arnold
\cite{A66}.

\subsection{The rotating rigid body}
The configuration space of a rigid body in $\R^3$ rotating around
its center of mass is the Lie group $SO(3)$.\footnote{See
\cite{MR1999} for further details on the material of this
subsection.} The corresponding Lie algebra is $\mathfrak{so}(3)$,
the space of antisymmetric $3\times 3$-matrices, which can be
identified with $\R^3$ via the map
$$\hat{}\ :\R^3\to\mathfrak{so}(3),\quad x=(x_1,x_2,x_3)\mapsto\hat x=\left(%
\begin{array}{ccc}
  0 & -x_3 & x_2 \\
  x_3 & 0 & -x_1 \\
  -x_2 & x_1 & 0 \\
\end{array}%
\right).$$ Let $I:\mathfrak{so}(3) \to \mathfrak{so}(3)^*$ be the
inertia matrix of the body. A \emph{left}-invariant metric
$\ska{\cdot}{\cdot}$ on $SO(3)$ is defined by setting
$$\ska{a}{b}=a\cdot Ib, \qquad  a,b\in \R^3 \simeq \mathfrak{so}(3),$$
at the identity, and extending it to all of $SO(3)$ by left
invariance. The basic observation is that $R(t)$ is a geodesic on
$(SO(3),\ska{\cdot}{\cdot})$ if and only if
$\hat\Omega(t)\coloneqq R(t)^{-1}\dot R(t)$ solves the classical
Euler equation for the motion of a rotating rigid body,
$$I\dot{\Omega} = (I\Omega)\times\Omega.$$
Physically, $\hat\Omega(t)$ represents the angular velocity in a
frame of reference fixed with respect to the body. The angular
velocity in the spatially fixed frame is given by $\dot
R(t)R(t)^{-1}$. In other words: Applying left and right
translations to the material angular velocity $\dot R(t)$, one
obtains the \emph{body} and the \emph{spatial} angular velocities,
which are both elements of the Lie algebra $\mathfrak{so}(3)$. The
body and spatial angular momenta, which are elements of the dual
$\mathfrak{so}(3)^*$, are given by $\Pi(t)=I\Omega(t)$ and
$\pi(t)=R(t)\Pi(t)$, respectively. The body and spatial quantities
are related by  the adjoint and coadjoint actions
\begin{equation}\label{rigidbodyadjoints}
  \hat\omega(t)= \Ad_{R(t)}\hat\Omega(t)=R(t)\hat\Omega(t)R(t)^{-1},\qquad\Pi(t)=\Ad^*_{R(t)}\pi(t).
\end{equation}
Conservation of (spatial) angular momentum implies that $\pi$ is in fact constant in
time, i.e.
\begin{equation}\label{piconserved}
  \frac{d\pi}{dt} = 0.
\end{equation}

\subsection{The CH equation}
For the CH equation
\begin{equation}\label{CHapp}
  u_t-u_{txx}+3uu_x=2u_xu_{xx}+uu_{xxx},\qquad x\in\S,\quad t \in \R,
\end{equation}
the configuration space is $G = \Diff(\S)$ with multiplication
$(\varphi,\psi)\mapsto\varphi\circ\psi$. Elements of the Lie
algebra $\mathfrak{g}$ are identified with functions $S^1 \to \R$.
A \emph{right}-invariant metric is defined by setting
$$\ska{u}{v}_{H^1}=\int_{\S} uAvdx=\int_{\S}(uv+u_xv_x)dx,$$
where $A=1-\partial_x^2:\mathfrak{g} \to \mathfrak{g}^*$ is the
inertia operator. The basic observation is that $\varphi(t)$ is a
geodesic in $(\Diff(\S),\ska{\cdot}{\cdot}_{H^1})$ if and only if
$u(t)=TR_{\varphi(t)^{-1}}\varphi_t(t)=\varphi_t(t)\circ\varphi(t)^{-1}$
satisfies (\ref{CHapp}). In other words, the CH equation is the
Euler equation on $(\Diff(\S),\ska{\cdot}{\cdot}_{H^1})$. Letting
$U = TL_{\varphi^{-1}}\varphi_t = (u\circ\varphi)\varphi_x^{-1}$,
$U$ and $u$ are the analogs of the body and spatial angular
velocities: they are obtained by left and right translation, respectively, of
the material velocity $\varphi_t$ to the Lie algebra. The momentum
in the spatial frame is $m=Au$. The analog of equation
(\ref{rigidbodyadjoints}) is
$$u(t) = \Ad_{\varphi(t)} U(t), \qquad m_0(t) = \Ad_{\varphi(t)}^*m(t),$$
where $m_0 = (m\circ\varphi)\varphi_x^2$ is the momentum in the
body frame. Since the metric now is right-invariant instead of
left-invariant, the analog of the conservation law
(\ref{piconserved}) is that the momentum $m_0$ in the body frame
is conserved,
$$\frac{dm_0}{dt} = 0, \qquad \text{i.e.} \qquad (m \circ \varphi)\varphi_x^2 = m_0.$$

\subsection{The 2CH equation}
For the 2CH equation (\ref{2CH}) the configuration space is the
semidirect product $G = \Diff(S^1) \circledS \mathcal{F}(S^1)$
introduced in Section \ref{Secsemidirect}. The Lie algebra
$\mathfrak{g}$ is identified with $\mathcal{F}(S^1) \times
\mathcal{F}(S^1)$. The inertia operator is $\text{diag}(A,\id)$
and the metric is the right-invariant metric $\langle \cdot, \cdot
\rangle$ defined in (\ref{metricdef}). The basic observation is
that $(\varphi(t), f(t))$ is a geodesic in $(\Diff(S^1) \circledS
\mathcal{F}(S^1),\ska{\cdot}{\cdot})$ if and only if
$$(u(t), \rho(t)) = TR_{(\varphi(t), f(t))^{-1}}(\varphi_t(t), f_t(t))$$
satisfies (\ref{2CH}). The analog of the body angular velocity is $(U_1, U_2) = TL_{(\varphi, f)^{-1}}(\varphi_t, f_t)$.
The spatial momentum is $(m, \rho) = (Au,\rho)$.
The analog of equation (\ref{rigidbodyadjoints}) is
$$(u(t), \rho(t)) = \Ad_{(\varphi(t), f(t))} (U_1(t), U_2(t))$$
and
$$(m_0(t), \rho_0(t)) = \Ad_{(\varphi(t),f(t))}^*(m(t), \rho(t))$$
where $(m_0, \rho_0)$ is the momentum in the body frame. In order
to find an explicit expression for $(m_0, \rho_0)$, we need to
compute the adjoint and coadjoint actions.

The adjoint action of $G$ on $\mathfrak{g} := T_{(\id,0)}G \simeq \mathcal{F}(S^1) \times \mathcal{F}(S^1)$ is defined by
$$\Ad_{(\varphi, f)}(v, \tau) := T_{(\id, 0)}I_{(\varphi, f)}\cdot (v, \tau), \qquad (v, \tau) \in \mathfrak{g},$$
where $I_{(\varphi, f)}:G \to G$ denotes the inner automorphism defined by
$$I_{(\varphi, f)}(\psi, g) = (\varphi, f)(\psi, g)(\varphi, f)^{-1}.$$
%= (\varphi \circ \psi \circ \varphi^{-1}, -f\circ\varphi^{-1} + (g + f \psi) \circ \varphi^{-1}).$$
A direct computation yields
\begin{align*}
 \Ad_{(\varphi,f)}(v,\tau)=(\Ad_\varphi v, (f_{x}v + \tau)\circ\varphi^{-1}), \qquad (v, \tau) \in \mathfrak{g},
\end{align*}
where $\Ad_\varphi v = (\varphi_xv)\circ\varphi^{-1}$ is the
adjoint action with respect to $\Diff(S^1)$. The $L^2$-pairing is
used to identify the (regular part of the) dual $\mathfrak{g}^*$
of $\mathfrak{g}$ with $\mathcal{F}(S^1) \times \mathcal{F}(S^1)$.
Since
\begin{align*}
\ska{(m, \rho)}{\Ad_{(\varphi,f)}(v, \tau)}&=\int_{\S}m\Ad_\varphi v dx+\int_{\S}\rho[(f_xv+\tau)\circ\varphi^{-1}]dx
    \\
&=\ska{\begin{pmatrix}
  (m\circ\varphi)\varphi_x^2+(\rho\circ\varphi)f_x\varphi_x \\
  (\rho\circ\varphi)\varphi_x
  \end{pmatrix}}{\begin{pmatrix}
  v \\
  \tau \end{pmatrix}},
\end{align*}
we find
\begin{align*}
& \Ad_{(\varphi,f)}^*(m, \rho)=\begin{pmatrix}
  (m\circ\varphi)\varphi_x^2+(\rho\circ\varphi)f_x\varphi_x \\
  (\rho\circ\varphi)\varphi_x
  \end{pmatrix}, \qquad (m,\rho) \in \mathfrak{g}^*.
\end{align*}
The analog of the conservation law (\ref{piconserved}) is that the momentum $(m_0, \rho_0)$ in the body frame is conserved,
$$\frac{d}{dt} \begin{pmatrix} m_0 \\ \rho_0 \end{pmatrix} = 0, \qquad \text{i.e.} \qquad
 \begin{pmatrix}
  (m\circ\varphi)\varphi_x^2+(\rho\circ\varphi)f_x\varphi_x \\
  (\rho\circ\varphi)\varphi_x
  \end{pmatrix} = \begin{pmatrix} m_0 \\ \rho_0
  \end{pmatrix}.$$
This explains the origin of the conservation law (\ref{ddtrhovarphi}) which was used in the proof of Proposition \ref{2CHlocalflowprop}.

\begin{center}
\begin{scriptsize}
\begin{tabular}{|l|c|c|c|}
  \hline
   & Rigid body & CH & 2CH \\
  \hline
  configuration space &$SO(3)$&$\Diff(\S)$&$\Diff(\S)\circledS \mathcal{F}(\S)$\\
  material velocity & $\dot R$ & $\varphi_t$ & $(\varphi_t,f_t)$\\
  spatial velocity & $\hat\omega=\dot RR^{-1}$ & $u = \varphi_t \circ \varphi^{-1}$ & $(u,\rho) = (\varphi_t\circ \varphi^{-1}, f_t\circ \varphi^{-1}) $ \\
  body velocity & $\hat\Omega=R^{-1}\dot R$ & $U=\frac{\varphi_t}{\varphi_x}$ & $\begin{pmatrix}
  U_1 \\
  U_2
\end{pmatrix}
  =\begin{pmatrix}
  \frac{\varphi_t}{\varphi_x} \\
  f_t - \frac{f_x}{\varphi_x}\varphi_t \end{pmatrix}$ \\
  inertia operator & $I$ & $A=1-\partial_x^2$ & $\begin{pmatrix}
  A & 0 \\
  0 & \text{id} \end{pmatrix}$ \\
  spatial momentum & $\pi=R\Pi$ & $m=Au$ & $(m, \rho) =(Au,\rho)$ \\
  body momentum & $\Pi=I\Omega$ & $m_0=(m\circ\varphi)\varphi_x^2$ & $\begin{pmatrix} m_0 \\ \rho_0 \end{pmatrix}=
  \begin{pmatrix}
  (m\circ\varphi)\varphi_x^2+(\rho\circ\varphi)f_x\varphi_x \\
  (\rho\circ\varphi)\varphi_x \end{pmatrix}$  \\
  spatial velocity (Ad) & $\hat\omega = \Ad_{R}\hat\Omega$ & $u = \Ad_{\varphi} U$ & $(u, \rho) = \Ad_{(\varphi, f)} (U_1, U_2)$\\
  body momentum (Ad*) & $\Pi=\Ad^*_{R}\pi$  & $m_0 =\Ad^*_\varphi m$ & $(m_0, \rho_0) =\Ad_{(\varphi,f)}^*(m, \rho)$ \\
  momentum conservation & $\pi = \text{const.}$ & $m_0 = \text{const.}$ & $(m_0, \rho_0) = \text{const.}$ \\
  \hline
\end{tabular}
\end{scriptsize}
\end{center}

\end{document}